\documentclass[12pt]{amsart}
\usepackage{amsmath}
\usepackage{amssymb}

\makeatletter

\def\O{{\mathcal O}}

\def\h{{\mathfrak h}}
\def\g{{\mathfrak g}}

\def\der{{\rm Der}\,}
\def\Bb#1{{\mathbb #1}}
\def\cal#1{{\mathcal #1}}

\newtheorem{Theorem}{Theorem}
\newtheorem{Lemma}[Theorem]{Lemma}

\newtheorem{Proposition}[Theorem]{Proposition}
\newtheorem{Definition}[Theorem]{Definition}

\theoremstyle{remark}
\newtheorem{Remark}[Theorem]{Remark}
\newtheorem{Remarks}[Theorem]{Remarks}
\newtheorem{Example}[Theorem]{Example}
\newtheorem{Question}[Theorem]{Question}

\begin{document}

\title[Euler characteristic of line singularities]
{A formula for Euler characteristic of line singularities on singular spaces}

\author{Guangfeng  Jiang}
\thanks{This article comes from a piece of the author's thesis which was
finished at Utrecht University advised by Prof. D. Siersma, to whom this author
is indebted very much for many discussions, questions and remarks.
This work was partially supported by JSPS, NNSFC, STCLN}
\address{{\it Present (till March 31, 2000)}: Department of  Mathematics,
    Faculty of Science, Tokyo Metropolitan University,   Minami-Ohsawa 1-1,
    Hachioji-shi,   Tokyo, 1920397 JAPAN  }
\email{jiang@comp.metro-u.ac.jp}

\address{{\it Permanent } Department of  Mathematics,
        Jinzhou Normal University,
        Jinzhou City,  Liaoning \ 121000, People's Republic of China  }
\email{jzjgf@mail.jzptt.ln.cn}


\begin{abstract}
We prove an algebraic formula for the Euler
characteristic of the Milnor fibres of functions with
critical locus a smooth curve on a space which is a weighted homogeneous
complete intersection with isolated singularity. 
\end{abstract}

\maketitle

\section{Introduction}

 For an analytic function germ $f: (X, 0)\longrightarrow ({\mathbb C},0) $
with critical locus $\Sigma \subset X$, there is a local Milnor
fibration induced by $f$. We are interested in the topology of the
Milnor  fibre $F$ of $f$ in the case when $\dim \Sigma =1$. It is well
known that in this case the homotopy type of $F$ is not necessarily a
bouquet of spheres in the middle dimension. 
The calculation of the Euler characteristic $\chi (F)$ of $F$ is of importance.
There  is a  so called Iomdin-L\^e
formula  \cite{Le 4}
which expresses the Euler characteristic of the Milnor fibre of $f$
by that of the {\it  series}  of $f$ with isolated singularities.

The question we are interested in 
 is that if there is a way to express $\chi (F)$ by some
``computable'' invariants determined only by   $(f,\Sigma ,X)$.
 When $X$ is $ {\mathbb C}^{m}$, the singular locus $\Sigma$
 of $f$ is a one dimensional complete intersection with 
 isolated singularity, and the transversal
 singularity type of $f$ along $\Sigma$ is Morse, 
 Pellikaan \cite{P4} answered this question positively. Pellikaan's
 formula expresses the Euler characteristic in terms of the
 Jacobian number $j(f)$, $\delta$ and the Milnor number $\mu(\Sigma)$
 of $\Sigma$.
These numbers can be computed directly by counting the dimensions of certain
finite dimensional vector spaces. 
 The development of computer algebra 
makes this kind of algebraic formulae  more and more important and popular.
 In this article we answer the question by proving a
 similar formula for function germs
 with line singularities  
on a  weighted homogeneous space $X$
with isolated complete intersection singularity
(see Proposition~\ref{maintheorem}).
Remark that for a function germ $f$ with
 isolated singularity on  a 
weighted homogeneous complete intersection with isolated singularity,
 Bruce and Robert \cite{BR} have  proved an
 algebraic formula for the Milnor number of $f$.

\section{Non-isolated singularities on singular spaces}

 Let ${\cal  O}_{{\mathbb C}^{m}}$ be
 the structure 
sheaf of ${\mathbb C}^{m}$. The stalk ${\cal  O}_{{\mathbb C}^{m},0}$
 of ${\cal  O}_{{\mathbb C}^{m}}$ at 0 is a local ring, consisting of
 germs at 0 of analytic functions on ${\mathbb C}^{m}$. 
The ring $ {\cal  O}_{{\mathbb    C}^{m},0}$ is often denoted by
 ${\cal  O}_{m}$,  or simply by ${\cal  O}$ when no confusion
 can be caused. The  unique maximal ideal of ${\cal  O}_m$ is  denoted  by 
 ${\mathfrak   m}_m$ or ${\mathfrak m}$.

Let $(X,0)\subset ({\mathbb C}^{m},0)$ be the germ of a reduced 
analytic subspace  $X$ of ${\mathbb    C}^{m}$  defined by an ideal 
${\mathfrak h}$ of ${\cal  O}$,
 generated by  $h_1, \ldots , h_p \in {\cal  O}$.
Let ${\mathfrak g} $
be the  ideal generated by $g_1, \ldots , g_n \in {\cal O}$. The germ  of the
 analytic space  defined by ${\mathfrak g}$ at $0$ is denoted by $(\Sigma,0)$.
 Write ${\cal O}_X:={\cal O}/{\mathfrak h}$ and  ${\cal O}_{\Sigma}:
={\cal O}/{\mathfrak g}.$

 Let $\der({\cal O})$ denote the ${\cal O}$-module
 of germs of analytic vector fields 
on ${\mathbb    C}^{m}$ at $0$. Then $\der({\cal O})$
 is a free ${\cal O}$-module with 
$\frac{\partial }{\partial z_1}, \ldots, \frac{\partial }{\partial z_m} $
as basis, where $z_1, \ldots , z_m$ are the local coordinates of
 $({\mathbb C}^{m}, 0).$ $\der({\cal O})$ is a Lie algebra with the
 bracket defined by
$[\xi, \eta]:=\xi\eta-\eta\xi$ for all $\xi, \eta \in \text{Der}.$

Define $\der _{\mathfrak h}({\cal O}):
=\{\xi\in \der ({\cal O})\mid \xi({\mathfrak h})\subset {\mathfrak h}\}$,
which is  the ${\cal  O}$-module of
 {\it logarithmic vector fields }
 along $(X,0)$  and  a Lie subalgebra of $\der({\cal O})$ \cite{S, BR}.
 When ${\mathfrak h}$ is a
radical ideal  defining the analytic space $X$, 
$\der _{\mathfrak h}({\cal O})$
 is often denoted by $D_X$.
 Geometrically, $D_X$ consists of  all the 
 germs of vector fields that  are tangent to the smooth part of $X$.
 When $X$ is a weighted homogeneous complete intersection with isolated 
singularity, one  can write down precisely all the generators of $D_X$
(cf. \cite{W}).
 
For $f\in {\cal  O}$, the ideal
 $J_X(f):=\left\{\xi(f)\mid \xi \in D_X \right\}$ 
is called the (relative) {\it Jacobian ideal } of $f$.
Obviously, when $X$ is the whole space ${\mathbb C}^{m}$, namely, 
${\mathfrak h}=\{0\}, $  then $J_X(f)=J(f)$, the Jacobian ideal of $f$.

Let ${\cal  S}=\{ S_{\alpha}\} $ be an analytic stratification 
of $X$, $f:(X, 0)\longrightarrow ({\mathbb    C},0)$ an analytic function germ.
 The {\it critical locus }
 $\Sigma ^{\cal  S}_f$
 of $f$ relative to the stratification ${\cal  S}$ is the union of  
the critical loci of $f$ restricted to each of the strata $S_{\alpha }$, 
namely, 
$\Sigma ^{\cal  S}_f
=\bigcup _{S_{\alpha}\in {\cal S} }\overline{\Sigma _{f|S_{\alpha }}}$.
We denote  $\Sigma ^{\cal  S}_f$ by $\Sigma _f$ when ${\cal  S}$ is clear 
from the context.
 If the dimension of $\Sigma ^{\cal  S}_f$ is not positive, we say 
that $f$  defines (or has)  
{\it isolated singularities} on $X$. If the dimension of 
$\Sigma ^{\cal  S}_f$  is 
positive,  we say that  $f$ defines (or has)
{\it non-isolated singularities} on $X$.
 If $\Sigma ^{\cal  S}_f$ is one 
dimensional smooth complex manifold, we say that  $f$ defines (or has)
 a {\it line singularity }  on $X$.

For  an analytic space $X$ embedded in a neighborhood $U$ of 
$0\in {\mathbb C}^m$,  there is  a {\it  logarithmic stratification } 
${\cal S}_{\log}:=\{S_{\alpha} \}$ of
 $U$   (see \cite{S,BR}).  
In general, ${\cal S}_{\log}$ is not   locally finite. 
If ${\cal S}_{\log}$ is locally finite, then   $X$ is said to be
 {\it holonomic}. 

Let $X$ be of pure dimension. The collection 
${\cal S}_{\log}'=\{X\cap S_{\alpha} \mid S_{\alpha} 
\in {\cal S}_{\log}  \}$ is a stratification of $X$ which  will be 
called the {\it logarithmic stratification } of $X$ in this article.
 Especially, when $X$ 
 has isolated singularity in $0$, then $\{0\}$ and the connected 
components of $X\setminus \{0\}$ form 
 a holonomic logarithmic stratification of $X$. So  0 is always a 
critical point of 
any germ $f:(X, 0)\longrightarrow ({\mathbb    C},0)$ relative to 
this stratification.  Hence for   $f\in {\mathfrak m}$, 
$\Sigma_f=\{p\in X\mid \xi (f)(p)=0 \text{ for all } \xi \in D_X\}$
 is the critical locus of $f$ relative to the logarithmic stratification.
 Obviously 
$\Sigma_f =X\cap {\cal  V}(J_X(f))$.

\begin{Definition}\label{primitiveideal}\rm
 Let ${\mathfrak h}=(h_1, \ldots, h_p)\subset 
{\mathfrak g}=(g_1, \ldots , g_n)$ be ideals of 
$ {\cal  O}_{{\mathbb    C}^m,0}$.
 Define a subset of ${\cal  O}$, called the 
{\it  primitive ideal} 
of ${\mathfrak g}$ relative to ${\mathfrak h}$:
$$\int _{\mathfrak h}{\mathfrak g}:=\{ f\in \g\mid \xi(f)
\subset {\mathfrak g}  \text{ for all } \xi \in {\rm Der}_{\h}({\cal O}) \}.$$
\end{Definition}

In the following we always assume that ${\mathfrak h}$ and ${\mathfrak g}$ are
 radical, $X={\mathcal V}(\h)$ and $ \Sigma ={\mathcal V}(\g).$ In this case,
 $\int _{\mathfrak h}{\mathfrak g}$ is denoted by $\int _X{\mathfrak g}$
or $\int {\mathfrak g}$ when no confusion can be caused by this.

\begin{Remarks}\label{remark1}\rm

\begin{itemize}
 \item[(1)] When $X$ is smooth  this  definition was given by Pellikaan
 \cite{P1,P2}.  It is straightaway to verify that 
$\int _{\mathfrak h}{\mathfrak g}$ is an ideal of ${\cal  O}$, and
${\mathfrak g}^2+{\mathfrak h}
\subset \int _{\mathfrak h}{\mathfrak g}\subset {\mathfrak g}$
 always holds. And for 
${\mathfrak g}_i\supset {\mathfrak h}$ $(i=1,2)$, we have 
$\int _{\mathfrak h}{\mathfrak g}_1 \cap \int _{\mathfrak h}{\mathfrak g}_2=
\int _{\mathfrak h}\left({\mathfrak g}_1\cap {\mathfrak g}_2\right)$;

\item[(2)] Geometrically, the relative primitive ideal collects all the
functions whose zero level surfaces pass through $\Sigma$ and are tangent
to the regular part $X_{\rm reg}$ of $X$ along 
$\Sigma \cap X_{\rm reg}$.

\item[(3)]  The singular locus (relative to ${\cal S}'_{\log}$) of  $f$ is
 $\Sigma _f:={\cal V}(J_X(f))\cap X$, and $\Sigma _f\subset f^{-1}(0)$ 
if $f(0)=0$. If	$f\in \int {\mathfrak g}$, then 
$\Sigma \subset X\cap {\cal  V}(J_X(f))=\Sigma _f$. Conversely,
for $f\in {\mathfrak m}$, we have $f\in \int {\mathfrak g}$ when 
$\Sigma \subset\Sigma _f$ and ${\mathfrak g}$ is radical. 
 The reason is: $J_X(f)\subset {\mathfrak g}$, and since $f$ takes finite
values on $\Sigma _f$ and $0\in \Sigma $, $f|_{\Sigma}=0$, so 
$f^k\in {\mathfrak g}$ for some $k\in {\mathbb    N}$. Hence 
$f\in {\mathfrak g}$ since
${\mathfrak g}$ is radical.

\item[(4)] The relative primitive ideals have been generalized to higher
 relative primitive ideals in \cite{JSimis}. Under the assumption that
 ${\mathfrak h}$ is pure dimension, ${\mathfrak g}$ is radical, and the
Jacobian ideal of ${\mathfrak h}$ is not contained in any associated prime of
${\mathfrak g}$, it was proved that the  primitive ideal
$\int _{\mathfrak h}{\mathfrak g}$ is the inverse image in ${\cal O}$
  of the second symbolic power of the quotient ideal
$\bar{\mathfrak g}:={\mathfrak g}/{\mathfrak h}$ of ${\cal O}_X$. 
Remark that the results in \cite{JSimis} generalized the results of 
\cite{Seibt, P1, P2}.
\end{itemize}
\end{Remarks}

\section{Transversal singularities}\label{transversal}

 Let $(X, 0)\subset ({\mathbb C}^{n+1},0)$ be the germ of a reduced analytic space  
with isolated singularity in 0.
 Let $\Sigma$ be a reduced curve germ on $(X,0)$ defined
by ${\mathfrak g}$ and have isolated singularity in 0. A germ 
$f\in \int {\mathfrak g}$
is called a 
{\it transversal $A_1$ singularity}
 along $\Sigma$ on $X$ 
if its singular locus $\Sigma _f= \Sigma$,  and, for $P\in \Sigma \setminus 0, f$
 has only $A_1$ singularity  transversal to the branch	of 
 $\Sigma$ containing $P$.
It was proved in \cite{J}, that
$f$ is a transversal $A_1$ singularity along $\Sigma$ on $X$
 if and only if the {\it Jacobian number }
 $j(f):=\dim(\g/(\h+J_X(f))< \infty$.

There exist  {\it  admissible linear forms }
 $l$ (see \cite{Le 4}) such that $\{l=0\}\cap
\Sigma =\{0\}$ and $\{l=0\}$ intersects  both $X$ 
 and ${ \Sigma}$ transversally, and $\{l=0\}\cap X\cap f^{-1}(0)$
 has isolated singularity at the origin.
We assume  $\{l=0\}$ is the first coordinate hyperplane $\{z_0=0\}$ of
 ${\mathbb C} ^{n+1}$. 

 Let $D_X$ be  generated by 
$\xi^0, \xi^1, \ldots \xi ^s$. Denote by
$D^0_X$ the submodule of $D_X$ generated by those $\xi ^i$
such that if we write
$\xi ^i =\sum_{j=0}^{n}\xi ^i_j{\frac{\partial } {\partial z_j}}$, then 
$ \xi ^i _0\notin {\mathfrak   g}$, and by
$D^1_X$ the submodule of $D_X$ generated by those $\xi ^i$
such that if we write
$\xi ^i =\sum_{j=0}^{n}\xi ^i_j{\frac{\partial } {\partial z_j}}$, then 
$ \xi _0\in {\mathfrak   g}$,
thus $D_X=D^0_X+D^1_X$. Denote
$$J^0(f):=D_X^0(f)=\{\xi(f)\mid \xi \in D_X^0\},\quad 
J^1(f):=D_X^1(f)=\{\xi(f)\mid \xi \in D_X^1\}.$$
 If $f$ is clear from the context we just write $J^0$ and $J^1$.

\begin{Lemma}\label{lemma1}\sl
 Let $f\in \int {\mathfrak   g}$, and $z_0, z_1,
 \ldots, z_n$  be the coordinates of ${\mathbb C} ^{n+1}$ such that $z_0=0$ 
 is admissible. The transversal singularity type of $f$ along every branch
of $\Sigma $ is constant at all the points of 
$\Sigma \setminus \{0\}$ if and only if 
$$\dim_{\mathbb C}\left(\frac{\cal O}{{\mathfrak g}
+((J^1+{\mathfrak h}):J^0)}\right)<\infty.$$
\end{Lemma}

\begin{proof} The inequality means that
 $\Sigma \cap {\cal  V}((J^1+{\mathfrak   h}):J^0)=\{0\}$.
For $\epsilon >0$ small enough, let $P\in \Sigma \cap \{z_0=t\}$, 
$0<|t|<\epsilon$. $P\notin  {\cal  V}((J^1+{\mathfrak   h}):J^0)$ 
if and only if 
$(J^0)_P\subset (J^1)_P$, where $(J^i)_P$ is the localization of the
 ideal at $P$. Since $z_0=0$ intersects both $\Sigma $ and $X$
transversally and $X$ is smooth at $P$, we can choose the local coordinates
 such that locally the branch of $\Sigma $ containing $P$
 is the first axis.
  Furthermore, we can
arrange the coordinate transformation such that under this transformation
$D_X^0$ and $D_X^1$ are preserved. By this we mean that
 any derivation of $X$ with the first
component non-zero at $P$ will remain non-zero at $P$ and any derivation
 of $X$ with first
component zero at $P$ will remain zero at $P$. We let $x=z_0, y_1, \ldots, 
y_{n-p}$ be the new local coordinates in a neighborhood of $P$ in $X$,
 then at $P$, 
$(J^0)_P=\left({\frac{\partial f}{\partial x}}\right){\cal  O}_{X,P}$ and
$(J^1)_P=\left({\frac{\partial f}{\partial y_1}},\ldots, 
{\frac{\partial f}{\partial y_{n-p}}}\right){\cal  O}_{X,P}$.
By \cite{P4}, the inequality is equivalent to the constancy of 
  the transversal singularity type of $f$.
\end{proof}

 Let $X$, $\Sigma $ and $f\in \int \g$ be the same as before.
 There is an integer
$k_0$ such that for all $k\geq k_0$, $f_k=f+{\frac{1}{k+1}}x^{k+1}$ 
defines an isolated  singularity at $O$. Let $\mu (f_k)$ be the Milnor
 number of $f_k$.  If $X$ is a weighted homogeneous complete intersection
 with isolated singularity, 
 by \cite{BR} 7.7, we have
$$\mu (f_k)=\dim \left({\frac{{\cal  O}}
{{\mathfrak   h}+J_X(f_k)}}
\right). $$
 
By the exact sequence
$$0\longrightarrow {\frac{{\mathfrak   g}+J_X(f_k)}{{\mathfrak   h}+J_X(f_k)}}
\longrightarrow {\frac{\cal  O}{{\mathfrak   h}+J_X(f_k)}}
\longrightarrow {\frac{\cal  O}{{\mathfrak   g}+J_X(f_k)}}
\longrightarrow 0,$$
we know that 
$$\mu (f_k)=\dim
 \left({\frac{{\mathfrak   g}+J_X(f_k)}{{\mathfrak   h}+J_X(f_k)}}\right)
+\dim \left({\frac{\cal  O}{{\mathfrak   g}+J_X(f_k)}}\right)
\eqno{(\ref{transversal}.1)}$$
Define 
$$e_k:=\dim {\frac{\cal  O}{{\mathfrak   g}+J_X(f_k)}}, \quad
\sigma(\Sigma X, 0):=\dim \left( {\frac{{\cal  O}_{\Sigma}}{J_X(x)}}\right),
\quad
 \text{multi}_x (\Sigma):
=\dim\left( {\frac{{\cal  O}_{\Sigma}}{(x)}}\right).
$$
Then obviously
$e_k=\sigma(\Sigma X, 0)+k \text{multi}_x (\Sigma), $
and 
$${\frac{{\mathfrak   g}+J_X(f_k)}{{\mathfrak   h}+J_X(f_k)}}=
{\frac{\mathfrak   g}{({\mathfrak   h}+J_X(f_k))\cap{\mathfrak   g}}}.$$
Hence
$$\mu (f_k) =  \sigma(\Sigma X, 0)+k\text{multi}_x (\Sigma) +
\dim\left( {\frac{\mathfrak   g}{({\mathfrak h}+J_X(f_k))\cap{\mathfrak   g}}}
\right)\eqno{(\ref{transversal}.2)}$$

\section{Line singularities}\label{linesing}

 In this  section, we assume that $\Sigma $ is a line in 
${\mathbb C} ^{n+1}$ defined by the ideal ${\mathfrak g}$, and
 $X$ is a space with isolated complete intersection singularity 
 defined by ${\mathfrak h}\subset {\mathfrak g}$.

\begin{Lemma}\label{decomposition} \sl
Let  $X $ be a space with isolated complete intersection singularity,
  containing a line
 $\Sigma$, which is chosen to be the first axis of a local coordinate
 system of  $\Bb C ^{n+1}$. If $X$ is weighted homogeneous with respect to 
 this coordinate system, then 
$D_X={\cal  O}\xi _E+D_X^1$ and 
\begin{align}
J_X(f_k)&=(\xi _E(f_k)){\cal  O}+J^1(f_k)\notag\cr
&\subset (\xi _E(f_k)){\cal  O}+ D_X^1(f)+D_X^1(x^{k+1}),\notag
\end{align}
where $J^1(f_k)=\{\xi(f_k)\mid \xi \in D_X^1\}\subset {\mathfrak   g}$.
\end{Lemma}
\begin{proof} Let $x, y_1, \ldots , y_n$
be  the local coordinates of $(\Bb C ^{n+1},0)$ in the statement of the lemma.
We know by \cite{W} that the Euler 
derivation $\xi _E=w_0x \frac{\partial}{\partial x}
+\sum\limits _{k=1}^{n}w_ky_k\frac{\partial}{\partial y_k}  \in D_X$, where 
$w_0, w_1, \ldots , w_n$ are the weights of the coordinates
$x, y_1, \ldots , y_n$ respectively.
Let $\xi =\xi ^0\frac{\partial}{\partial x}+
 \sum\limits _{k=1}^{n}\xi ^k\frac{\partial}{\partial y_k}\in D_X$. 
 Since 
${\cal  O}_{\Sigma} $ is a principal ideal domain and $X$ has isolated
 singularity at $O$, $\bar{\xi} ^0\in (\bar x)$.
 Let $\bar{\xi}^0=\bar{x}\bar{\xi}_1^0$.
Then $\xi '=\xi -\frac{1}{w_0}\xi _1^0\xi _E$ has
the coefficient of
${\frac{\partial }{\partial x}}$ in ${\mathfrak   g}$. 
\end{proof}

\begin{Lemma}\label{basic}\sl 
Let $\Sigma$ and $X$ be the same as in Lemma~\ref{decomposition}.
 Let $f\in \int {\mathfrak   g}$ have 
transversal $A_1$
singularity along  $\Sigma$. For $ k>0$ sufficiently large, we have
$$\tilde{J}:=\xi _E(f){\mathfrak   g}+D_X^1(f)+{\mathfrak   h}=
({\mathfrak   h}+J_X(f_k))\cap{\mathfrak   g}$$
\end{Lemma}
\begin{proof}
 By Lemma~\ref{decomposition}, we have
$$
({\mathfrak   h}+J_X(f_k))\cap{\mathfrak   g}
=(\xi _E(f_k))\cap {\mathfrak g}+D_X^1(f_k)+{\mathfrak h}.
$$
For $a\in (\xi _E(f_k))\cap {\mathfrak   g}$, 
$a=a_0x^{k+1}+a_0\xi _E(f)\in {\mathfrak   g}$.
 Since $\xi _E(f)\in {\mathfrak   g}$,
 $a_0\in {\mathfrak   g}$. Hence 
$ (\xi _E(f_k))\cap {\mathfrak   g}=(\xi _E(f_k)){\mathfrak   g}$ and
$$\begin{array}{ll}
({\mathfrak   h}+J_X(f_k))\cap{\mathfrak   g}
&=(x^{k+1}+\xi _E(f)) {\mathfrak g}+D_X^1(f_k)+{\mathfrak h}\cr
&\subset (x^{k+1}+\xi _E(f)) {\mathfrak   g}+x^kD^1_X(x)
 +D^1_X(f)+{\mathfrak   h}
\end{array}
\eqno{(\ref{linesing}.1)}$$

Since $j(f) <\infty $,
 there is a $k_1$ such that when $k>k_1$, $x^k{\mathfrak   g}\subset 
{\mathfrak   h}+J_X(f)
=(\xi _E(f))+D^1_X(f)+{\mathfrak   h}$.
 
By Lemma~\ref{lemma1}, there is a $k_2$ such that 
$x^{k_2}\in {\mathfrak   g}+((D^1_X(f)+{\mathfrak   h}):(\xi _E (f)))$, and
$x^{k_2}\xi _E (f)\in (\xi _E(f)){\mathfrak   g}+D^1_X(f)+{\mathfrak   h}$.
Hence when $k>k_1+k_2$, 
$$a_0(x^{k+1}+\xi _E(f))+b_0x^k\in \tilde J\quad \text{ and }\quad
({\mathfrak   h}+J_X(f_k))\cap{\mathfrak   g}\subset \tilde J.$$

 On the other hand, there is an integer $n_1>>0$  such that when $k\geq n_1$
(see (\ref{transversal}.2))
$$x^k{\mathfrak   g}\subset ({\mathfrak   h}+J_X(f_k))\cap{\mathfrak g}
\subset \tilde{J}\subset 
{\mathfrak   g}.$$
Then for $k>n_1$,  by the equality in (\ref{linesing}.1) 
\begin{align}
\tilde{J}&
\subset  \xi _E\left(f+{\frac{x^{k+1}}{k+1}}\right)
{\mathfrak g}+x^{k+1}{\mathfrak g}
+D_X^1(f)+{\mathfrak   h}\notag\cr
&= (\xi _E(f)+x^{k+1}){\mathfrak g}+D_X^1(f)+{\mathfrak h}
+x^{k+1}{\mathfrak g}\notag\cr
&=({\mathfrak h}+J_X(f_k))\cap{\mathfrak g}+x^{k+1}{\mathfrak g}\notag\cr
&\subset({\mathfrak h}+J_X(f_k))\cap{\mathfrak g}
+x^{k}{\mathfrak g}\notag\cr  
&\subset({\mathfrak h}+J_X(f_k))\cap{\mathfrak g}+x^{k-n_1}\tilde{J}\notag\cr
&\subset({\mathfrak h}+J_X(f_k))\cap{\mathfrak g}+{\mathfrak m}\tilde{J}.\notag
\end{align}
By Nakayama's lemma,  
$({\mathfrak   h}+J_X(f_k))\cap{\mathfrak   g}=\tilde{J}$.
\end{proof}

\begin{Lemma}\label{reduction}\sl
 Under the assumption of Lemma~\ref{basic}, we have
 \begin{xxalignat}{2}
1)&\quad  L:={\frac{{\mathfrak   h}+J_X(f)}{\tilde{J}}}\cong
 {\frac{(\xi _E(f))}{(\xi _E(f))\cap\tilde{J}}};\cr
2)&\quad {\rm Ann}(L)={\mathfrak g}+((D_X^1(f)+{\mathfrak h}):(\xi _E(f)));\cr
3) &\quad L\cong  \frac{\cal  O}{{\rm Ann}(L)}
={\frac{\cal  O}{\g+((D_X^1(f)+{\mathfrak   h}):(\xi _E(f)))}}.
\end{xxalignat}
\end{Lemma}
\begin{proof} It is an easy exercise in commutative algebra,
we omit it.
\end{proof}
 
 From the exact sequence
$$0\longrightarrow {\frac{{\mathfrak   h}+J_X(f)}{({\mathfrak   h}
+J_X(f_k))\cap{\mathfrak   g}}}
\longrightarrow{\frac{\mathfrak   g}{({\mathfrak   h}
+J_X(f_k))\cap{\mathfrak   g}}}
\longrightarrow{\frac{\mathfrak   g}{{\mathfrak   h}+J_X(f)}}
\longrightarrow 0,$$
we have
$$\mu (f_k) =j(f) +e_k+
\dim\left({\frac{{\mathfrak   h}+J_X(f)}{({\mathfrak   h}
+J_X(f_k))\cap{\mathfrak   g}}} \right).$$ 

Notice  that  $\sigma(\Sigma X, 0)=\text{multi}_x(\Sigma)=1$, $e_k=k+1$. By
Lemma~\ref{reduction}, we have
$$\mu (f_k)=k+1+j(f) 
+\dim \left(\frac{\cal  O}{\g+ ((D_X^1(f)+{\mathfrak   h}):
(\xi _E(f)))} \right)$$

Let $F$ and $F_k$ are  the Milnor fibre of
$f$ and $f_k$ respectively. Iomdin-L\^e's formula  \cite{Le 4, Ti},
 says that 
$$\chi (F)=\chi (F_k)+(-1)^{\dim X}(k+1).$$
But
$\chi (F_k)=1+(-1)^{\dim X-1}\mu (f_k)$.
Since in our case $\sigma(\Sigma X, 0)=1$,  we have proved

\begin{Proposition}\label{maintheorem}\sl
Let  $X $ be a space with isolated  complete intersection 
 singularity  containing a line $\Sigma$, which is chosen to be the first
 axis of a local coordinate system of  $\Bb C ^{n+1}$. Assume that
 $X$ is weighted homogeneous with respect to 
 the coordinate system.
 For an analytic function $f\in \int _X\g$ with $j(f)<\infty$, 
 the Euler characteristic of the Milnor fibre $F$ of $f$ is
$$
\chi (F)=1+(-1)^{\dim X-1}(j(f) +\nu)\eqno{(\ref{linesing}.2)}$$
 where 
$$\nu =\dim \left( {\frac{\cal  O}{\g+((D_X^1(f)+{\mathfrak   h}):
(\xi _E(f)))}} \right).\eqno{\square}$$
\end{Proposition}

\begin{Remark}\rm Formula (\ref{linesing}.2) allows us to use a
 computer program  to compute 
the Euler characteristic effectively. In fact, we have a small  {\it
Singular }\cite{Singular} program to calculate $\chi(F)$ for a function
with critical locus  a line   on a hypersurface $X$. We use it to check the 
 examples in Example~\ref{ex}.
\end{Remark}

 Let $\Sigma $ be a line in 
${\mathbb C} ^{n+1}$ defined by the ideal ${\mathfrak g}=(y_1, \ldots , y_n)$,
 $X$  a space with isolated complete intersection singularity
 defined by ${\mathfrak h}=(h_1,\ldots, h_p)\subset {\mathfrak g}$.
By changing the generators of $\g$, we can write 
$$h_i\equiv \sum_{k=1}^p b_{ik}y_k\mod \g ^2
$$
such that the determinant $b$ of the matrix $B=(b_{ik})$ is a 
non-zero divisor in ${\cal  O}_{\Sigma }$.
Note that $y_1, \ldots y_p$ are
projected to zero or the generators of the torsion part $T(M)$  of
 the conormal module 
$M={\mathfrak   g}/({\mathfrak g}^2+{\mathfrak h})$
of $\g/\h$,
and $y_{p+1}, \ldots y_n$ form a basis of the free module $N=M/T(M)$.
   We call 
$\lambda(\Sigma X):=\dim_{\Bb C} T(M)=\dim \O/(\g+(b))$
the {\it  torsion number} of the space pair $(\Sigma, X)$.
For  $f=\sum h_{kl}y_ky_l\in  {\mathfrak   g}^2$,  define
$$\Delta =\det(h_{kl})_{p+1\le k,l \le n},\quad
\delta _f:=\dim(\frac{\O}{\g+\Delta}).$$

\begin{Question}\label{question}\rm
 For $f\in \g ^2$, does the equality
 $\nu=2\lambda(\Sigma X)+\delta _f-1$
always hold? Or under what conditions does it hold?
\end{Question}

The equality in  Question~\ref{question} and (\ref{linesing}.2)
 show the geometric meaning of $\chi(F)$.
 By using  {\it Singular }\cite{Singular} we have checked that it holds
 for all the examples we know.

\begin{Example}\label{ex}\rm
 Let $X_{k,l}$ be  an  $A_{k,l}$ singularity defined by
 $h=x^ly+x^{s}z^2+yz$. This is an $A_k$ singularity with $k=2l+s-1$
($l\ge 1, s\ge 0$). Since we take the line $\Sigma$ with torsion number $l$
as $x$-axis, we have the definition equation (see \cite{JS}). 
By resolution of singularities, one can prove that for 
any function $f$ on $X_{k,l}$, if it has  isolated line singularity and
 $A_1$ type transversal singularity, then the Milnor fibre $F$ of $f$
 is a bouquet of circles (see [11]). Then $\mu (f)=j(f)+\nu$.

\begin{itemize}
\item[(1)]
We consider a function $g: X _{k,l}\longrightarrow {\mathbb C}.$
For generic $(a,b,c)\in \Bb C^3$, let $g=ay^2+byz+cz^2$, for example,
take $g=y^2-yz+\frac{1}{2}z^2$, a calculation shows that $\mu(g)=6l-3$.
\item[(2)] Let
$f:X _{k,l}\longrightarrow {\mathbb C} $ be defined by $f=y+\frac{1}{2}z^2$.
 A calculation shows that
 $$\mu(f)=\left\{
\begin{array}{ll}
l-1 & \text{ if } s=0, \cr 
l+3s-2 & \text{ if } 1\le s\le  l-1, \cr 
4l-2 & \text{ if } s\ge  l.
\end{array}\right.$$
\end{itemize}
\end{Example}

 \end{document}